\documentclass[12pt,reqno,a4paper]{amsart} 
\usepackage[breaklinks]{hyperref}
\usepackage{mathtools} 

\usepackage[headings]{fullpage} 

\usepackage{pslatex} 

\usepackage{microtype} 
 
\renewcommand{\qedsymbol}{$\square$}
\newenvironment{Proof}[1][Proof]{\par\noindent\textbf{#1.}~}
{\hfill\qedsymbol\smallskip\par}

\newcommand{\dx}{\mathrm{d}}
\newcommand{\eps}{\varepsilon}

\newcommand{\Stilde}{\widetilde{S}}
\newcommand{\Etilde}{\widetilde{E}}

\newcommand{\second}{\prime\prime}

\newcommand{\Odip}[2]{\mathcal{O}_{#1}\!\left(#2\right)\mathchoice{\!}{}{}{}}
\newcommand{\Odipg}[2]{\mathcal{O}_{#1}\Bigl(#2\Bigr)\mathchoice{\!}{}{}{}}
\newcommand{\Odipm}[2]{\mathcal{O}_{#1}\bigl(#2\bigr)\mathchoice{\!}{}{}{}}
\newcommand{\Odig}[1]{\mathcal{O}\Bigl(#1\Bigr)\mathchoice{\!}{}{}{}}

\newcommand{\Odim}[1]{\mathcal{O}\bigl(#1\bigr)}
\newcommand{\Odi}[1]{\Odip{}{#1}}
\newcommand{\odip}[2]{{o}_{#1}\!\left(#2\right)\mathchoice{\!}{}{}{}}
\newcommand{\odi}[1]{\odip{}{#1}}
 
\makeatletter
\newenvironment{Biggcases}{%
  \matrix@check\Biggcases\env@Biggcases
}{%
  \endarray %
}
\def\env@Biggcases{%
  \let\@ifnextchar\new@ifnextchar
  \Biggl\lbrace
  \def\arraystretch{1.2}%
  \array{@{}l@{\quad}l@{}}%
}
\makeatother

\newtheoremstyle{sltheorems}
{10pt}
{6pt}
{\slshape}
{}
{\bfseries}
{.}
{.5em}
{\thmname{#1}\thmnumber{ #2}\thmnote{ (#3)}}

\theoremstyle{sltheorems} 
\newtheorem{Theorem}{Theorem}
\newtheorem{Lemma}{Lemma}

\allowdisplaybreaks
\begin{document}
\title[asymptotic formulae for binary problems, I: density $3/2$]{Short intervals asymptotic formulae for binary problems with primes and powers, I: density $3/2$} 
\author[]{Alessandro Languasco \lowercase{and} Alessandro Zaccagnini}

\begin{abstract}
We prove that suitable asymptotic formulae in short intervals hold for the problems
of representing an integer as a sum of a prime and a square, or a prime square.
Such results are obtained both assuming the Riemann Hypothesis and
in the unconditional case.
\end{abstract}

\maketitle

\allowdisplaybreaks  
\section{Introduction}

In this first paper devoted to study asymptotic formulae in short intervals for additive problems with primes and squares, we focus our attention on density-$3/2$ problems, \emph{i.e.}, on representing integers as sum of a prime and a square. 
In the forthcoming paper \cite{LanguascoZ2015d} we will consider density-$1$ problems.

Let $\eps>0$, $N$ be a sufficiently large integer and  let further $H$ be an integer such that $N^{\eps}<H=\odi{N}$ as $N \to \infty$.
Taking $n\in [N, N+H]$, the key quantities  are
\[ 
r^{\prime}_{1,2}(n) = 
 \sum_{p+m^2=n} \log p
 \quad
\textrm{and}
\quad
r^{\second}_{1,2}(n) = 
 \sum_{p_{1}+p_{2}^2=n} 
\log p_{1} \log p_{2}.
\]
Since it is well known that the expected behaviour of such functions is erratic,  
to work in a more regular situation  we will study their average asymptotics  over a suitable short interval.
We write $f=\infty(g)$ for $g=\odi{f}$.
In the following we prove
\begin{Theorem}
\label{asymp-part1}  Assume the Riemann Hypothesis (RH) holds.
Then
\[
\sum_{n = N+1}^{N + H} r^{\prime}_{1,2}(n)
=
H N^{1/2} 
+ 
 \Odim{N^{3/4}(\log N)^2+H^{3/2}(\log N)^{3/2}+ HN^{1/3}\log N}
\]
as $N\to \infty$ uniformly for   $\infty(N^{1/4}(\log N)^{2})\le H \le \odi{N/(\log N)^{3}}$.
\end{Theorem}

\begin{Theorem}
\label{asymp-unc-part1}
Let $\eps>0$. Then there exist two constants $C=C(\eps)>0$, 
$C_{1}=C_{1}(\eps)>0$ such that  
\[
\sum_{n = N+1}^{N + H}  r^{\prime}_{1,2}(n)
=
H N^{1/2}
+
 \Odig{(H^{1/2} N^{3/4} +   H N ^{1/2}) \exp \Big( -C \Big( \frac{\log N}{\log \log N} \Big)^{1/3} \Big)}
\]
as $N\to \infty$ uniformly for 
\[
 N^{1/2} \exp \Big( - C_{1} \Big( \frac{\log N}{\log \log N} \Big)^{1/3}\Big) \le H \le N^{1-\eps}.
\]
\end{Theorem}

A direct trial following the lines of  Lemma 11 of Plaksin \cite{Plaksin1981}
leads to have a square summand in $[N,N+H]$ and hence the final 
uniformity range  has to be larger than $H>N^{1/2}$ 
which is weaker than our results above. 

Concerning the sum of a prime and a prime square we have the following
\begin{Theorem}
\label{asymp-part2}  
Assume the Riemann Hypothesis holds. Then
\[
\sum_{n = N+1}^{N + H} r^{\second}_{1,2}(n)
=
H N^{1/2} 
+ 
 \Odig{\frac{H^2}{N^{1/2}}+  N^{3/4}(\log N)^3+ HN^{1/3}(\log N)^2}
\]
as $N\to \infty$ uniformly for $\infty(N^{1/4}(\log N)^{3}) \le H \le \odi{N}$. 
\end{Theorem}

\begin{Theorem}
\label{asymp-unc-part2}
Let $\eps>0$. Then there exists a constant $C=C(\eps)>0$ such that  
\[
\sum_{n = N+1}^{N + H}  r^{\second}_{1,2}(n)
= 
H N^{1/2}
 +
 \Odig{H N ^{1/2} \exp \Big( -C \Big( \frac{\log N}{\log \log N} \Big)^{1/3} \Big)}
\]
as $N\to \infty$ uniformly for $N^{7/12+\eps}\le H \le N^{1-\eps}$.
\end{Theorem}
 
In this case too, 
a direct trial following the lines of  Lemma 11 of  Plaksin \cite{Plaksin1981}
leads to weaker uniformity ranges:  
$H\gg N^{3/4} (\log N)^A$ assuming RH and  $H\gg N^{7/24+1/2+\eps}$ unconditionally.

Our results are proved via a circle method technique; in fact
for Theorem \ref{asymp-part2} we'll need  the original Hardy-Littlewood approach 
(using infinite series instead of finite sums) otherwise Lemma \ref{App-BCP-Gallagher}
below requires $H> N^{1/2}$. 
This is similar to the phenomenon we already encountered in our paper
\cite{LanguascoZ2016a}. We also remark that the original 
Hardy-Littlewood approach can be applied in proving Theorem \ref{asymp-part1}
too; but in this case it will just lead to replace the error term
$H^{3/2}(\log N)^{3/2}$ with the slightly better one  $H^2 N^{-1/2}$.

Clearly our result implies the existence of an integer represented
as a sum of a prime and a square, or a prime square, in the stated intervals.
Concerning this we have to remark that  Kumchev and Liu  \cite{KumchevL2009} unconditionally
proved the existence of an integer which is the sum of a prime and a prime  square
in the shorter interval $H>N^{0.33}$ but without any information about the
relevant asymptotic formula. As far as we know this is the best known result
for the the sum of a prime and a square case too.
  
\medskip
{\bf Acknowledgements.}    This research was partially supported by the grant PRIN2010-11 \textsl{Arithmetic Algebraic Geometry and Number Theory}. We wish to thank the referee for his/her remarks.

\section{Definitions and Lemmas}

Let $L=\log N$,  $r_{0}(m)$ be the number of representations of $m$ as a sum of two squares (recall that $r_{0}(m)\ll m^{\eps}$ is a well-known fact) and
\[
R^{\prime}_{1,2}(n) =
\sum_{\substack {m_{1}+m_{2}^2=n\\1\leq m_{1}, \,  m_{2}^2\leq N}}
\Lambda(m_{1})
\quad
\textrm{and}
\quad
R^{\second}_{1,2}(n) =
\sum_{\substack {m_{1}+m_{2}^2=n\\1\leq m_{1}, \, m_{2}^2\leq N}}
\Lambda(m_{1})\Lambda(m_{2}).
\]

As $n\in [N,N+H]$, $N\to \infty$ and $H=\odi{N}$, it is easy to see that
\begin{align}
\notag
r^{\prime}_{1,2}(n)  
&=
 \sum_{\substack {p+m^2=n\\1\leq p,m^{2}\leq N}} \log p
 +\Odi{\frac{HL}{N^{1/2}}+H^{1/2}L}
 =
R^{\prime}_{1,2}(n)  +
\Odig{
\sum_{\substack {p^{j}+m^{2}=n\\1\leq p^{j},\, m^{2}\leq N;\, j\ge 2}} \log p
}
 +\Odim{H^{1/2}L}
\\
\notag
&
=
R^{\prime}_{1,2}(n) +
\Odig{
\sum_{\substack {p^{2k}+m^{2}=n\\1\leq p^{2k},\, m^{2}\leq N;\, k\ge 1}} \log p
+
\sum_{\substack {p^{2k+1}+m^{2}=n\\1\leq p^{2k+1},\, m^{2}\leq N;\, k\ge 1}} \log p
}
 +\Odim{H^{1/2}L}
\\
\label{r-R-reduction1}
&
=
R^{\prime}_{1,2}(n)+
\Odi{r_0(n)L^{2}  +  n^{1/3}L + H^{1/2}L}
=
R^{\prime}_{1,2}(n)  +\Odim{n^{1/3}L+ H^{1/2}L},
\end{align}
using the Prime Number Theorem,
and, similarly, that
\begin{equation}
\label{r-R-reduction2}
r^{\second}_{1,2}(n)  
=
 R^{\second}_{1,2}(n)+\Odim{n^{1/3}L^{2}+ H^{1/2}L^{2}}.
\end{equation}

So from now on we can work with the uppercase-$R$ functions.
 Let  now $\ell\ge 1$ be an integer and 
\begin{align}
\notag
S_{\ell}(\alpha) = \sum_{1\leq m^{\ell} \leq N}\Lambda(m) e(m^{\ell} \alpha), & \quad 
T_{\ell} (\alpha) = \sum_{m^{\ell} \le N} e(m^{\ell} \alpha) , \\
\label{main-defs}
f_{2}(\alpha) =\frac{1}{2}\sum_{1 \leq m\leq N} m^{-1/2}e(m\alpha), 
& \quad
U(\alpha,H) = \sum_{1\leq m\leq H}e(m\alpha),
\end{align}
where  $e(\alpha) = e^{2\pi i\alpha}$.
We also have the usual numerically explicit inequality
\begin{equation}
\label{UH-estim}
\vert U(\alpha,H) \vert
\le
\min \bigl(H;  \vert \alpha\vert ^{-1}\bigr),
\end{equation}
see, \emph{e.g.}, on page 39 of Montgomery \cite{Montgomery1994}.
Let further 
\begin{equation}
\label{B-def}
B=B(N,c)= \exp \Big( c   \Big( \frac{L}{\log L} \Big)^{1/3} \Big),
\end{equation} 
where $c=c(\eps)>0$ will be chosen later.

In the proofs we will need the following lemmas.
In fact we will use them just for $\ell=1,2$ but we take this occasion
to describe the  general case. 
We explicitly remark that for $\ell=1$ the  proof of 
Lemma \ref{zac-lemma} gives just trivial results; in this case
a non-trivial estimate, which, in any case, is not useful in this context, 
can be obtained following the line of Corollary 3 of \cite{LanguascoP1994}.
\begin{Lemma} 
\label{zac-lemma}
Let $\ell\ge 2$ be an integer and $0<\xi\leq 1/2$. Then
\[
\int_{-\xi}^{\xi} \vert T_{\ell}(\alpha) \vert ^2\, \dx \alpha 
=
2 \xi N^{1/\ell} + 
\begin{Biggcases}
\Odi{L} & \text{if}\ \ell =2\\
\Odip{\ell}{1} & \text{if}\ \ell > 2
\end{Biggcases}
\]
and
\[
\int_{-\xi}^{\xi} \vert S_{\ell}(\alpha) \vert ^2\, \dx \alpha 
=
\frac{2\xi}{\ell} N^{1/\ell} L +   \Odipm{\ell}{\xi N^{1/\ell}} + 
\begin{Biggcases}
\Odim{L^{2}} & \text{if}\ \ell =2\\
\Odip{\ell}{1} & \text{if}\ \ell > 2.
\end{Biggcases}
\]
\end{Lemma}
\begin{Proof}
By symmetry we can  integrate over $[0,\xi]$.
We use Corollary 2 of Montgomery-Vaughan  \cite{MontgomeryV1974}  
with $T=\xi$, $a_r=1$ and $\lambda_r= 2\pi r^\ell$ thus getting
\begin{align*}
\int_{0}^{\xi} \vert T_{\ell}(\alpha) \vert ^2\, \dx \alpha 
&=
\sum_{r^{\ell} \le N} \bigl(\xi +\Odim{\delta_r^{-1}}\bigr)
=
\xi N^{1/\ell} +  \Odi{\xi}
+\Odipg{\ell}{ 
\sum_{r^{\ell} \le N} r^{1-\ell}} 
\end{align*}
since $\delta_{r} = \lambda_r - \lambda_{r-1} \gg_{\ell} r^{\ell-1}$. 
The last error term is $\ll_{\ell}1$ if $\ell >2$ and $\ll L$ otherwise.
This proves the first part of Lemma \ref{zac-lemma}.
Arguing analogously with  $a_r=\Lambda(r)$, by the Prime Number Theorem we get
\begin{align*}
\int_{0}^{\xi} \vert S_{\ell}(\alpha) \vert ^2\, \dx \alpha 
&=
\sum_{r^{\ell} \le N} \Lambda(r)^2 \bigl(\xi +\Odim{\delta_r^{-1}}\bigr)
=
\frac{\xi}{\ell} N^{1/\ell} L  +  \Odipm{\ell}{\xi  N^{1/\ell}}
+\Odipg{\ell}{ 
\sum_{r^{\ell} \le N} \Lambda(r)^2 r^{1-\ell}}\, .
\end{align*}
Again by the Prime Number Theorem, the last error term is $\ll_{\ell}1$ if $\ell >2$ and $\ll L^{2}$ otherwise.
The second part of Lemma \ref{zac-lemma} follows.
\end{Proof}

We need the following lemma which collects the results 
of Theorems 3.1-3.2 of \cite{LanguascoZ2013b}; see also Lemma 1 of \cite{LanguascoZ2016a}.
\begin{Lemma}  
\label{App-BCP-Gallagher}
Let $\ell  > 0$ be a real number and $\eps$ be an arbitrarily small
positive constant. Then there exists a positive constant 
$c_1 = c_{1}(\eps)$, which does not depend on $\ell$, such that
\[
\int_{-1/K}^{1/K}
\vert
S_{\ell}(\alpha) - T_{\ell}(\alpha)  
\vert^2 
\, \dx \alpha
\ll_{\ell}   N^{2/\ell -1}
\Bigl(
\exp \Big( - c_{1}  \Big( \frac{L}{\log L} \Big)^{1/3} \Big)
+
\frac{K L^{2}}{N}
\Bigr),
\]
uniformly for  $N^{1-5/(6\ell)+\eps}\le K \le N$.  
Assuming further RH we get 
\[
\int_{-1/K}^{1/K}
\vert
S_{\ell}(\alpha) - T_{\ell}(\alpha) 
\vert^2 
\, \dx \alpha
\ll_{\ell} 
\frac{N^{1/\ell} L^{2}}{K} + K N^{2/\ell-2} L^{2},
\]
uniformly for  $N^{1-1/\ell}\le K \le N$.  
\end{Lemma}

\section{Proof of Theorem \ref{asymp-part1}}  
{}From now on, we denote 
\(
E_{\ell}(\alpha) : =S_\ell(\alpha) - T_\ell(\alpha).
\)
By \eqref{main-defs} it is an easy matter to see that
\begin{align}
\notag
\sum_{n=1}^H &
R^{\prime}_{1,2}(n+N) = 
\int_{-1/2}^{1/2} S_{1}(\alpha)T_{2}(\alpha) U(-\alpha,H)e(-N\alpha) \, \dx \alpha
\\
\notag
 &=
 \int_{-1/2}^{1/2}\!\!\! T_{1}(\alpha) f_{2}(\alpha) U(-\alpha,H)
 e(-N\alpha)\, \dx \alpha 
 +
\int_{-1/2}^{1/2}\!\!\! T_{1}(\alpha) ( T_{2}(\alpha) - f_{2}(\alpha) ) U(-\alpha,H)
 e(-N\alpha)\, \dx \alpha
 \\
\label{approx-th1-part1} 
 &\hskip 1cm 
 + 
 \int_{-1/2}^{1/2}  E_{1}(\alpha) T_{2}(\alpha) U(-\alpha,H)
e(-N\alpha)\, \dx \alpha
= I_{1}+I_{2} + I_{3}, 
\end{align}
say.
Now we evaluate $I_{1}$. 
A direct calculation and Lemma 2.9 of Vaughan \cite{Vaughan1997} give
\begin{align}
\notag
\int_{-1/2}^{1/2} & T_{1}(\alpha) f_{2}(\alpha)    e(-(n+N)\alpha)\, \dx \alpha 
=
\frac{1}{2}\sum_{\substack {m_{1}+m_{2} =n+N\\ 1\leq m_{1}, m_{2} \leq N}}m_{1}^{-1/2}
=
\frac{1}{2}\sum_{1\le m \leq N} (n+N-m)^{-1/2}
\\
\label{M-eval-part1}
&=
\frac{\Gamma(1/2)}{2\Gamma(3/2)} (n+N)^{1/2} 
+\Odim{  n^{1/2}}
= 
 (n+N)^{1/2}  + \Odim{n^{1/2}}.
 \end{align} 
By \eqref{approx-th1-part1}-\eqref{M-eval-part1} we obtain
\begin{equation}
\label{I1-eval-part1} 
I_{1}
=  \sum_{n=1}^H (n+N)^{1/2} + \Odim{H^{3/2}} =
H N^{1/2} +  \Odim{H^{3/2}} . 
\end{equation}

Now we estimate $I_{2}$. 
We first recall,  by Theorem 4.1 of Vaughan \cite{Vaughan1997},  that
$\vert T_{2}(\alpha) -f_{2}(\alpha) \vert \ll (1+\vert \alpha \vert N)^{1/2}$. 
Using  also the inequality
$ T_{1}(\alpha)  \ll  \min(N; \vert \alpha \vert^{-1})$, we get
\begin{align}
\notag
I_{2}
&\ll
\int_{-1/2}^{1/2} \vert T_{1}(\alpha) \vert  \vert T_{2}(\alpha) -f_{2}(\alpha) \vert 
\vert U(\alpha,H) \vert \, \dx \alpha
\\
\label{I2-estim-part1}
& 
\ll
H N \int_{-1/N}^{1/N} \, \dx \alpha
+
H N^{1/2}\int_{1/N}^{1/H} \frac{\dx \alpha}{\alpha^{1/2}}  
+
N^{1/2}\int_{1/H}^{1/2} \frac{\dx \alpha}{\alpha^{3/2}}  
 \ll 
H^{1/2}N^{1/2}.
\end{align}

To estimate $I_{3}$ we need  Lemmas \ref{zac-lemma}-\ref{App-BCP-Gallagher}.
By \eqref{UH-estim} and the Cauchy-Schwarz inequality we have
\begin{equation*}
I_{3} 
\ll 
\Bigl(\int_{-1/2}^{1/2}  
\!\!\! \vert E_{1}(\alpha)  \vert^{2}  \min(H; \vert \alpha \vert^{-1}) \, \dx \alpha\Bigr)^{1/2} 
\Bigl(\int_{-1/2}^{1/2} \!\!\! \vert T_{2}(\alpha) \vert^{2}  \min(H; \vert \alpha \vert^{-1}) \, \dx \alpha \Bigr)^{1/2} 
 =
(J_{1} J_{2})^{1/2},
\end{equation*}
say.
Since
\begin{equation*}
\label{J1-split-part1} 
J_{1} 
\ll 
H\int_{-1/H}^{1/H}  \vert E_{1}(\alpha)  \vert^{2}\, \dx \alpha 
+
\int_{1/H}^{1/2}
 \vert E_{1}(\alpha)  \vert^{2}\frac{\, \dx \alpha}{ \alpha},
\end{equation*}
by Lemma \ref{App-BCP-Gallagher} with $\ell=1$  and partial integration we get 
\begin{equation}
\label{J1-estim-part1} 
J_{1} \ll  NL^{3} +  H^{2} L^{2}.
\end{equation}
Arguing analogously and using Lemma \ref{zac-lemma} with $\ell=2$, we obtain
\begin{equation}
\label{J2-estim-part1} 
J_{2} \ll  (N^{1/2}  + H) L.
\end{equation}
Hence combining  \eqref{J1-estim-part1}-\eqref{J2-estim-part1} we have
\begin{equation}
\label{I3-estim-part1} 
I_{3} 
\ll 
N^{3/4}L^{2}
+
H N^{1/4} L^{3/2}
+
H^{1/2} N^{1/2} L^{2}
+
H^{3/2}  L^{3/2} .
\end{equation}
Now using \eqref{approx-th1-part1}, \eqref{I1-eval-part1}-\eqref{I2-estim-part1} and \eqref{I3-estim-part1}, 
we  can finally write 
\begin{equation*} 
    \sum_{n=1}^H R^{\prime}_{1,2}(n+N) 
    =
    H N^{1/2} +  \Odim{N^{3/4}L^{2}+H^{3/2} L^{3/2}+ H^{1/2}N^{1/2}L^{2} + H N^{1/4} L^{3/2}}.
\end{equation*}
Using \eqref{r-R-reduction1},
Theorem \ref{asymp-part1}  hence follows for $\infty(N^{1/4}L^{2}) \le H \le \odi{N/L^{3}}$.
\qed

\section{Proof of Theorem \ref{asymp-unc-part1}}  

We need now to split the main interval in a different way. Recalling \eqref{B-def} and 
\(
E_{\ell}(\alpha) =S_\ell(\alpha) - T_\ell(\alpha),
\)
by \eqref{main-defs} we have
\begin{align} 
\notag
\sum_{n = N+1}^{N + H}&  R^{\prime}_{1,2}(n)  
=
\int_{-B/H}^{B/H} S_{1}(\alpha)T_{2}(\alpha)  U(-\alpha,H)  e(-N\alpha) \, \dx \alpha
+
\int\limits_{\mathclap{[-1/2,-B/H]\cup  [B/H,1/2]}} 
S_{1}(\alpha)T_{2}(\alpha)  U(-\alpha,H) e(-N\alpha) \, \dx \alpha
\\
\notag
 &\!\!\!\!\!\!=
 \int_{-B/H}^{B/H}\!\!\!\!  T_{1}(\alpha) f_{2}(\alpha)   U(-\alpha,H) 
 e(-N\alpha)\, \dx \alpha 
 +\!\!
 \int_{-B/H}^{B/H} \!\!\! \! T_{1}(\alpha) ( T_{2}(\alpha) - f_{2}(\alpha) )    U(-\alpha,H) 
 e(-N\alpha)\, \dx \alpha 
 \\
 \notag
 &\hskip1cm + 
 \int_{-B/H}^{B/H}  E_{1}(\alpha) T_{2}(\alpha)  U(-\alpha,H)  
e(-N\alpha)\, \dx \alpha
+
\int\limits_{\mathclap{[-1/2,-B/H]\cup  [B/H,1/2]}} 
S_{1}(\alpha)T_{2}(\alpha)  U(-\alpha,H) e(-N\alpha) \, \dx \alpha
\\
\label{approx-th2-part1} 
&
= I_{1}+I_{2} + I_{3}+ I_{4}, 
\end{align}
say.
Arguing as in \eqref{M-eval-part1}, using \eqref{UH-estim} and
$f_2(\alpha)\ll \min(N^{1/2}, 1/\vert \alpha\vert^{1/2})$
(see Lemma 2.8 of Vaughan \cite{Vaughan1997}), we obtain 
\begin{equation} 
        \label{I1-unc-eval-part1} 
I_{1}
=  
\sum_{n=1}^{H} (n+N)^{1/2} 
+ 
\Odim{H^{3/2}}
+ 
\Odig{\int_{B/H}^{1/2} \frac{\dx \alpha}{\alpha^{5/2}} }
=
H N^{1/2} +  \Odim{H^{3/2}} . 
\end{equation}

$I_{2}$ can be estimate as in \eqref{I2-estim-part1} and gives
\begin{equation} 
\label{I2-unc-estim-part1}
I_{2}
 \ll 
H^{1/2}N^{1/2}.
\end{equation} 

Now we estimate $I_{3}$. 
By \eqref{UH-estim} the Cauchy-Schwarz inequality we have
\[
I_{3}
\ll 
H
\Bigl(\int_{-B/H}^{B/H}  
 \vert E_{1}(\alpha) \vert^{2}  \, \dx \alpha\Bigr)^{1/2} 
\Bigl(\int_{-B/H}^{B/H}  \vert T_{2}(\alpha) \vert^{2}    \, \dx \alpha \Bigr)^{1/2} 
 =
H(J_{1} J_{2})^{1/2},
\]
say. 
By Lemma \ref{App-BCP-Gallagher}  we can write that 
\begin{equation}
\label{J1-estim-unc-part1} 
J_{1} \ll N 
\exp \Big( - c_{1}  \Big( \frac{L}{\log L} \Big)^{1/3} \Big)  
\end{equation}
provided that $ N^{-1-\eps/2}< B/H < N^{-1/6-\eps/2}$; hence
  $N^{1/6+\eps}\le H \le N^{1-\eps}$ suffices.  
By Lemma \ref{zac-lemma}  with $\ell=2$, we obtain
\begin{equation}
\label{J2-estim-unc-part1} 
J_{2} \ll  \frac{N^{1/2} B}{H}+ L.
\end{equation}

Hence combining  \eqref{J1-estim-unc-part1}-\eqref{J2-estim-unc-part1}
for   $N^{1/6+\eps}\le H \le N^{1-\eps}$  we have
\begin{equation}
\label{I3-estim-unc-part1} 
I_{3}
\ll 
(H^{1/2} N^{3/4} B^{1/2} + H N ^{1/2} L^{1/2})
\exp \Big( - \frac{c_{1}}{2}  \Big( \frac{L}{\log L} \Big)^{1/3} \Big).
\end{equation}

Now we estimate $I_{4}$. By \eqref{UH-estim}, the Prime Number Theorem, Lemma \ref{zac-lemma} with $\ell=2$
and a partial integration argument
we get
\begin{align}
\notag
I_{4}
&\ll
\int_{B/H}^{1/2}
\vert  S_1(\alpha) T_{2}(\alpha) \vert
\frac{\dx \alpha}{\alpha}
\ll
\Bigl(
\int_{B/H}^{1/2}  
\vert S_1(\alpha) \vert^2 \frac{\dx \alpha}{\alpha}
\Bigr)^{1/2}
\Bigl(
\int_{B/H}^{1/2}  
\vert T_{2}(\alpha)\vert^2\frac{\dx \alpha}{\alpha}
\Bigr)^{1/2}
\\
\label{I4-estim-unc-part1}
&\ll
\Bigl(
\frac{HNL}{B}
\Bigr)^{1/2}
\Bigl(
 N^{1/2} 
+
\frac{H L}{B} 
+
\int_{B/H}^{1/2}  
(\xi N^{1/2}  + L) \frac{\dx \xi}{\xi^2}
\Bigr)^{1/2} 
\ll
\frac{HN^{1/2} L}{B} + 
\frac{H^{1/2} N^{3/4} L}{B^{1/2}}  .
\end{align}

Now using \eqref{approx-th2-part1}-\eqref{I2-unc-estim-part1} and \eqref{I3-estim-unc-part1}-\eqref{I4-estim-unc-part1}  and choosing $0<c<c_{1}$ in \eqref{B-def},
we have that there exists a constant $C=C(\eps)>0$ such that
\begin{equation*} 
\sum_{n = N+1}^{N + H}   R^{\prime}_{1,2}(n) 
    =
   H N^{1/2} 
 +
 \Odig{(H^{1/2} N^{3/4} +   H N ^{1/2}) \exp \Big( -C \Big( \frac{L}{\log L} \Big)^{1/3} \Big)}
\end{equation*}
uniformly for for  $N^{1/6+\eps}\le H \le N^{1-\eps}$.
Using \eqref{r-R-reduction1},
Theorem \ref{asymp-unc-part1}  hence follows for 
\[
 N^{1/2} \exp \Big( - C_{1} \Big( \frac{L}{\log L} \Big)^{1/3}\Big) \le H \le N^{1-\eps}
\]
for every $0<C_{1}= C_{1}(\eps)<2C$.
\qed 
  
\section{Proof of Theorem \ref{asymp-part2} }  

We need  the original Hardy-Littlewood approach otherwise Lemma \ref{App-BCP-Gallagher} 
implies that we need to assume $H\ge N^{1/2}$. 
Let further 
\begin{align} 
\label{tilde-main-defs}
\Stilde_\ell(\alpha) = \sum_{n=1}^{\infty} \Lambda(n) e^{-n^{\ell}/N} e(n^{\ell}\alpha),
\
\widetilde{R}^{\second}_{1,2}(n) =
\sum_{m_{1}+m_{2}^2=n}
\Lambda(m_{1})\Lambda(m_{2})
\quad
\textrm{and}
\
z= 1/N-2\pi i\alpha.
\end{align} 
{}From now on, we denote 
\(
\Etilde_{\ell}(\alpha) : =\Stilde_\ell(\alpha) -  \Gamma(1/\ell)/(\ell z^{1/\ell}).
\)
We remark  
\begin{equation}
\label{z-estim}
\vert z\vert ^{-1} \ll \min \bigl(N, \vert \alpha \vert^{-1}\bigr)
\end{equation}
and, arguing analogously to  \eqref{r-R-reduction1}-\eqref{r-R-reduction2}, that
\begin{equation}
\label{r-R-reduction3}
r^{\second}_{1,2}(n) 
=
\widetilde{R}^{\second}_{1,2}(n) 
+ \Odim{n^{1/3}L^{2}}.
\end{equation}

By \eqref{tilde-main-defs} it is an easy matter to see that 
\begin{align}
\notag 
\sum_{n = N+1}^{N + H}&  e^{-n/N}
\widetilde{R}^{\second}_{1,2}(n) 
= 
\int_{-1/2}^{1/2} \Stilde_{1}(\alpha)\Stilde_{2}(\alpha) U(-\alpha,H)e(-N\alpha) \, \dx \alpha
\\
\notag
 &=
 \int_{-1/2}^{1/2} \frac{\pi^{1/2}}{2 z^{3/2}}U(-\alpha,H)
 e(-N\alpha)\, \dx \alpha
+
 \int_{-1/2}^{1/2}   \frac{1}{z} \Etilde_{2}(\alpha)  U(-\alpha,H)
 e(-N\alpha)\, \dx \alpha 
 \\
\notag
 &\hskip0.7cm + 
 \int_{-1/2}^{1/2} \frac{\pi^{1/2}}{2 z^{1/2}}  \Etilde_{1}(\alpha)  U(-\alpha,H)
e(-N\alpha)\, \dx \alpha
 + 
 \int_{-1/2}^{1/2}  \Etilde_{1}(\alpha)
\Etilde_{2}(\alpha)  U(-\alpha,H)
e(-N\alpha)\, \dx \alpha 
\\
&
\label{approx-th1-part2} 
 = I_{1}+I_{2} + I_{3} + I_{4}, 
\end{align}
say.
We evaluate $I_{1}$.  Using Lemma~4 of \cite{LanguascoZ2016a} 
we immediately get  
\begin{equation}
\label{I1-eval-part2}
I_{1}
=
 \frac{\pi^{1/2}}{2\Gamma(3/2)}    \sum_{n = N+1}^{N + H}  n^{1/2}e^{-n/N} + \Odi{\frac{H}{N}}
  =
\frac{H N^{1/2}}{e}  + \Odi{\frac{H^2}{N^{1/2}}}.
\end{equation}

Now we estimate $I_{2}$.  
By \eqref{z-estim}, the Cauchy-Schwarz inequality and Lemma 3 of \cite{LanguascoZ2016a}, we obtain
\begin{align}
\notag
I_{2}
&\ll 
HN
 \int_{-1/N}^{1/N}  \vert \Etilde_{2}(\alpha)  \vert    \, \dx \alpha 
 +
 H
 \int_{1/N}^{1/H}  \vert \Etilde_{2}(\alpha)  \vert   \, \frac{\dx \alpha}{\alpha}
 +
 \int_{1/H}^{1/2}  \vert \Etilde_{2}(\alpha)  \vert  \, \frac{\dx \alpha}{\alpha^{2}}
\\
\notag
& 
\ll
HN^{1/4}L 
+ 
H
\Bigl( \int_{1/N}^{1/H} \!\!\!\!  \vert \Etilde_{2}(\alpha)  \vert^2   \frac{\dx \alpha}{\alpha}\Bigr)^{1/2}\Bigl( \int_{1/N}^{1/H}  \frac{\dx \alpha}{\alpha}\Bigr)^{1/2}
+
\Bigl( \int_{1/H}^{1/2}\!\!\!\!  \vert \Etilde_{2}(\alpha)  \vert^2 \frac{\dx \alpha}{\alpha}\Bigr)^{1/2}\Bigl(  \int_{1/H}^{1/2}  \frac{\dx \alpha}{\alpha^{3}}\Bigr)^{1/2}
\\
\label{I2-estim-part2}
& 
\ll
HN^{1/4}L + HN^{1/4}L^{2} + HN^{1/4}L^{3/2}
\ll
HN^{1/4}L^{2} .
\end{align} 

Now we estimate $I_{3}$.  
By \eqref{z-estim}, the Cauchy-Schwarz inequality and Lemma 3 of \cite{LanguascoZ2016a}, we have
\begin{align}
\notag
I_{3}
&\ll 
HN^{1/2}
 \int_{-1/N}^{1/N}  \vert \Etilde_{1}(\alpha)  \vert    \, \dx \alpha 
 +
 H
 \int_{1/N}^{1/H}  \vert \Etilde_{1}(\alpha)  \vert   \, \frac{\dx \alpha}{ \alpha ^{1/2}}
 +
 \int_{1/H}^{1/2}  \vert \Etilde_{1}(\alpha)  \vert  \, \frac{\dx \alpha}{\alpha^{3/2}}
\\
\notag
& 
\ll
H L 
+ H 
\Bigl( \int_{1/N}^{1/H}  \vert \Etilde_{1}(\alpha)  \vert^2   \, \frac{\dx \alpha}{\alpha}\Bigr)^{1/2}\Bigl( \int_{1/N}^{1/H} \dx \alpha\Bigr)^{1/2}
+
\Bigl( \int_{1/H}^{1/2}  \vert \Etilde_{1}(\alpha)  \vert^2   \, \frac{\dx \alpha}{\alpha}\Bigr)^{1/2}\Bigl(  \int_{1/H}^{1/2}  \frac{\dx \alpha}{\alpha^{2}}\Bigr)^{1/2}
\\
\label{I3-estim-part2}
& 
\ll
H L + H^{1/2}N^{1/2}L^{3/2}  
\ll
H^{1/2}N^{1/2}L^{3/2}.
\end{align} 
By \eqref{UH-estim} and the Cauchy-Schwarz inequality  we can write
\begin{align*}
I_{4}
&\ll 
H
\Bigl(\int_{-1/H}^{1/H}  
 \vert \Etilde_{1}(\alpha)  \vert^{2}  \, \dx \alpha\Bigr)^{1/2} 
\Bigl(\int_{-1/H}^{1/H}  \vert \Etilde_{2}(\alpha) \vert^{2}    \, \dx \alpha \Bigr)^{1/2} 
\\
&\hskip1cm
+
\Bigl(\int_{1/H}^{1/2}  
 \vert \Etilde_{1}(\alpha) \vert^{2}   \, \frac{\dx \alpha}{\alpha}\Bigr)^{1/2} 
\Bigl(\int_{1/H}^{1/2}  \vert \Etilde_{2}(\alpha) \vert^{2}    \, \frac{\dx \alpha}{\alpha} \Bigr)^{1/2} 
=
J_{1} + J_{2},
\end{align*}
say.
By Lemma 3 of \cite{LanguascoZ2016a} and partial integration on $J_{2}$, we obtain 
\begin{equation*} 
J_{1} \ll  N^{3/4}L^{2}
\quad 
\textrm{and}
\quad
J_{2} \ll  N^{3/4}L^{3}
\end{equation*}
and hence we get 
\begin{equation}
\label{I4-estim-part2}
I_{4}
\ll  
N^{3/4}L^{3}.
\end{equation}
Now using \eqref{approx-th1-part2}-\eqref{I3-estim-part2}
and \eqref{I4-estim-part2} we  can  write 
\begin{equation}
\label{almost-done} 
   \sum_{n = N+1}^{N + H} e^{-n/N}
\widetilde{R}^{\second}_{1,2}(n)
    = 
    \frac{H N^{1/2}}{e} + 
 \Odig{\frac{H^2}{N^{1/2}}+H^{1/2}N^{1/2}L^{3/2} +  N^{3/4}L^{3}}
\end{equation}
which is an asymptotic relation for $\infty(N^{1/4}L^{3}) \le H \le \odi{N}$.
{}From \eqref{r-R-reduction3} and  $e^{-n/N}=e^{-1}+ \Odi{H/N}$ for 
$n\in[N+1,N+H]$, we get 
\begin{equation}
\label{th3-final} 
   \sum_{n = N+1}^{N + H} 
r^{\second}_{1,2}(n)
    = 
    H N^{1/2} 
    + 
 \Odig{\frac{H^2}{N^{1/2}}+   N^{3/4}L^{3} + H N^{1/3}L^{2}}
 +
 \Odig{\frac{H}{N}\sum_{n = N+1}^{N + H} 
   \widetilde{R}^{\second}_{1,2}(n)}.
\end{equation}
Using  $e^{n/N}\leq e^{2}$ 
and \eqref{almost-done} for $H$ in the previously mentioned range,
it is easy to see that the last error term is
$\ll H^{2}N^{-1/2}$.
Combining \eqref{th3-final} and the last remark,
Theorem \ref{asymp-part2}  hence follows for $\infty(N^{1/4}L^{3}) \le H \le \odi{N}$.
\qed

\section{Proof of Theorem \ref{asymp-unc-part2}}  

In the unconditional case  we can use the 
finite sums approach. Recalling \eqref{main-defs}-\eqref{B-def}
and 
\(
E_{\ell}(\alpha) =S_\ell(\alpha) - T_\ell(\alpha),
\)
we have
\begin{align} 
\notag
\sum_{n = N+1}^{N + H} & R^{\second}_{1,2}(n) 
=
\int_{-B/H}^{B/H} S_{1}(\alpha)S_{2}(\alpha)  U(-\alpha,H)  e(-N\alpha) \, \dx \alpha
+
\int\limits_{\mathclap{[-1/2,-B/H]\cup  [B/H,1/2]}} 
S_{1}(\alpha)S_{2}(\alpha)  U(-\alpha,H) e(-N\alpha) \, \dx \alpha
\\
\notag
 &=
 \int_{-B/H}^{B/H}T_{1}(\alpha) T_{2}(\alpha) U(-\alpha,H)
 e(-N\alpha)\, \dx \alpha
 +
 \int_{-B/H}^{B/H} S_{1}(\alpha) E_{2}(\alpha) U(-\alpha,H)
  e(-N\alpha)\, \dx \alpha 
 \\
 \notag
 &\hskip0.3cm + 
 \int_{-B/H}^{B/H}  E_{1}(\alpha)T_{2}(\alpha)  U(-\alpha,H)
 e(-N\alpha)\, \dx \alpha
+
\int\limits_{\mathclap{[-1/2,-B/H]\cup  [B/H,1/2]}} 
S_{1}(\alpha)S_{2}(\alpha)  U(-\alpha,H) e(-N\alpha) \, \dx \alpha
\\
\label{approx-th2-part2} 
&= I_{1}+I_{2} + I_{3}+ I_{4}, 
\end{align}
say.
Using
$\vert T_{2}(\alpha) - f_{2}(\alpha) \vert \ll (1+\vert \alpha \vert N)^{1/2}$
(by Theorem 4.1 of Vaughan \cite{Vaughan1997}) and
$ T_{1}(\alpha)  \ll  \min(N; \vert \alpha \vert^{-1})$
we obtain
\begin{align}
\notag
 I_{1}
 &=
  \int_{-B/H}^{B/H} \!\!\! T_{1}(\alpha) f_{2}(\alpha)  U(-\alpha,H)
 e(-N\alpha)\, \dx \alpha
 +
  \int_{-B/H}^{B/H} \!\!\! T_{1}(\alpha) (T_{2}(\alpha)-f_{2}(\alpha)) U(-\alpha,H)
   e(-N\alpha)\, \dx \alpha
 \\
 \label{I1-frakI1-confronto}
 &
 =
{\mathfrak I}_{1} + \Odig{
 H N \int_{-1/N}^{1/N}\!\!\!    \dx \alpha
+
H N^{1/2}\int_{1/N}^{1/H}\!\! \!   \frac{\dx \alpha}{\alpha^{1/2}}
+
N^{1/2}\int_{1/H}^{B/H}\! \!\!   \frac{\dx \alpha}{\alpha^{3/2}}
} 
=
{\mathfrak I}_{1} + \Odim{H^{1/2}N^{1/2}},
\end{align}
say.
Using \eqref{UH-estim} and arguing as in \eqref{M-eval-part1} we obtain
\begin{align}
\notag
{\mathfrak I}_{1} 
& =
\sum_{n=1}^{H} 
\int_{-1/2}^{1/2} T_{1}(\alpha) f_{2}(\alpha)    e(-(n+N)\alpha)\, \dx \alpha 
+
\Odig{\int_{B/H}^{1/2}   \frac{\dx \alpha}{\alpha^{5/2}}}
\\
\label{M-eval-part2}
&
= 
\sum_{n=1}^{H} (n+N)^{1/2}  
 + 
 \Odig{\sum_{n=1}^{H} 
n^{1/2}}
+
\Odig{\frac{H^{3/2}}{B^{3/2}}}.
\end{align}

By \eqref{I1-frakI1-confronto}-\eqref{M-eval-part2} we obtain
\begin{equation}
\label{I1-unc-eval-part2} 
I_{1}
=  
\sum_{n = N+1}^{N + H} n^{1/2} 
+ 
\Odim{ H^{3/2} + H^{1/2}N^{1/2}} 
=
 H N^{1/2} + \Odim{ H^{3/2} + H^{1/2}N^{1/2}} . 
\end{equation}
Now we estimate $I_{2}$. 
By the Cauchy-Schwarz inequality we can write
\begin{align*}
I_{2}
\ll
H
\Bigl(\int_{-B/H}^{B/H}
\vert E_{2}(\alpha)\vert^{2}  \, \dx \alpha\Bigr)^{1/2} 
\Bigl(\int_{-B/H}^{B/H} \vert S_{1}(\alpha) \vert^{2}   \, \dx \alpha \Bigr)^{1/2} 
=
H(J_{1} J_{2})^{1/2},
\end{align*} 
say.
By Lemma \ref{App-BCP-Gallagher}  we get 
\[
J_{1} \ll  \exp \Big( - c_{1}  \Big( \frac{L}{\log L} \Big)^{1/3} \Big),
\]
provided that $N^{-1-\eps/2}<B/H<N^{-7/12-\eps/2}$; hence   $N^{7/12+\eps}\le H \le N^{1-\eps}$
suffices.
By the Prime Number Theorem we obtain
$ J_{2} \ll  N L $ and hence
\begin{equation}
\label{I2-estim-part2-unc}
I_{2}
\ll
H N^{1/2} L^{1/2}\exp \Big( - \frac{c_{1}}{2}  \Big( \frac{L}{\log L} \Big)^{1/3} \Big)
\ll
H N^{1/2} \exp \Big( - \frac{c_{1}}{4}  \Big( \frac{L}{\log L} \Big)^{1/3} \Big),
\end{equation}
uniformly for  $N^{7/12+\eps}\le H \le N^{1-\eps}$.

Now we estimate $I_{3}$. 
By the Cauchy-Schwarz inequality we have
\[
I_{3}
\ll 
H
\Bigl(\int_{-B/H}^{B/H}
 \vert E_{1}(\alpha) \vert^{2} \, \dx \alpha\Bigr)^{1/2} \!
\Bigl(\int_{-B/H}^{B/H} \vert T_{2}(\alpha) \vert^{2}   \, \dx \alpha \Bigr)^{1/2} 
 =
 H(K_{1} K_{2})^{1/2},
\]
say.
By Lemma \ref{App-BCP-Gallagher}  we get 
\begin{equation}
\label{K1-estim-unc-part2} 
K_{1} \ll N 
\exp \Big( - c_{1}  \Big( \frac{L}{\log L} \Big)^{1/3} \Big)
\end{equation}
provided that $N^{-1-\eps/2}<B/H < N^{-1/6-\eps/2}$; hence  $N^{1/6+\eps}\le H \le N^{1-\eps}$ suffices.  
By Lemma \ref{zac-lemma}  with $\ell=2$, we obtain
\begin{equation}
\label{K2-estim-unc-part2} 
K_{2} \ll   \frac{N^{1/2} B}{H}  +  L.
\end{equation}
Hence combining  \eqref{K1-estim-unc-part2}-\eqref{K2-estim-unc-part2}
for $N^{1/6+\eps}\le H \le N^{1-\eps}$ we have
\begin{equation}
\label{I3-estim-unc-part2} 
I_{3}
\ll 
(H^{1/2} N^{3/4} B^{1/2}+ H  N^{1/2}L^{1/2})
\exp \Big( - \frac{c_{1}}{2}  \Big( \frac{L}{\log L} \Big)^{1/3} \Big).
\end{equation}

Now we estimate $I_{4}$. By \eqref{UH-estim}, the Prime Number Theorem, Lemma \ref{zac-lemma} with $\ell=2$
and a partial integration argument
we get 
\begin{align}
\notag
I_4
&\ll
\int_{B/H}^{1/2}
\vert  S_1(\alpha) S_{2}(\alpha) \vert
\frac{\dx \alpha}{\alpha}
\ll
\Bigl(
\int_{B/H}^{1/2}  
\vert S_1(\alpha) \vert^2 \frac{\dx \alpha}{\alpha}
\Bigr)^{1/2}
\Bigl(
\int_{B/H}^{1/2}  
\vert S_{2}(\alpha)\vert^2\frac{\dx \alpha}{\alpha}
\Bigr)^{1/2}
\\
\notag
&\ll
\Bigl(
\frac{HNL}{B}
\Bigr)^{1/2}
\Bigl(
 N^{1/2} L
+
\frac{HL^{2}}{B}
+
L
\int_{B/H}^{1/2}  
(\xi N^{1/2}  + L) \frac{\dx \xi}{\xi^2}
\Bigr)^{1/2}
\\
\label{I4-estim-unc-part2}
&\ll
\Bigl(
\frac{HNL}{B}
\Bigr)^{1/2} \Bigl(
 N^{1/2} L^{2}
+
\frac{HL^{2}}{B}
\Bigr)^{1/2} .
\end{align}

Now using \eqref{approx-th2-part2}, 
 \eqref{I1-unc-eval-part2}-\eqref{I2-estim-part2-unc} and \eqref{I3-estim-unc-part2}-\eqref{I4-estim-unc-part2}, 
and choosing $0<c<c_{1}$ in \eqref{B-def},
we have that there exists a constant $C=C(\eps)>0$ such that
\begin{equation*} 
\sum_{n = N+1}^{N + H} R^{\second}_{1,2}(n) 
    =
   H N^{1/2} 
 +
 \Odig{ H N ^{1/2} \exp \Big( -C \Big( \frac{L}{\log L} \Big)^{1/3} \Big)}
\end{equation*}
uniformly for $N^{7/12+\eps}\le H \le N^{1-\eps}$.
Using \eqref{r-R-reduction2},
Theorem \ref{asymp-unc-part2}  hence follows for $N^{7/12+\eps}\le H \le N^{1-\eps}$.
\qed

\renewcommand{\bibliofont}{\normalsize}

\bigskip
\noindent
Alessandro Languasco, Dipartimento di Matematica, Universit\`a
di Padova, Via Trieste 63, 35121 Padova, Italy. {\it e-mail}: languasco@math.unipd.it

\medskip
\noindent
Alessandro Zaccagnini, Dipartimento di Matematica e Informatica, Universit\`a di Parma, Parco
Area delle Scienze 53/a, 43124 Parma, Italy. {\it e-mail}: alessandro.zaccagnini@unipr.it
\end{document}